\documentclass{amsart}
\usepackage[dvips]{graphicx}
\usepackage{times}

\newtheorem{theorem}{Theorem}[section]
\newtheorem{lemma}[theorem]{Lemma}

\newtheorem{corollary}[theorem]{Corollary}

\theoremstyle{definition}
\newtheorem{definition}[theorem]{Definition}
\newtheorem{remark}[theorem]{Remark}

\numberwithin{equation}{section}

\newcommand{\zz}{\mathbb{Z}}

\newcommand{\cpk}{\mathbb{CP}^2}
\newcommand{\cpkk}{\overline{\mathbb{CP}}{}^2}

\begin{document}

\title{Symplectic Tori in Rational Elliptic Surfaces}

\author{Tolga Etg\"u}
\address{Department of Mathematics and Statistics, McMaster University,
Hamilton, Ontario, L8S 4K1, Canada}
\email{etgut@math.mcmaster.ca}

\author{B. Doug Park}
\address{Department of Pure Mathematics, University of Waterloo, Waterloo,
Ontario, N2L 3G1, Canada}
\email{bdpark@math.uwaterloo.ca}
\thanks{B.D. Park was partially supported by NSERC and CFI grants.}

\date{July 28, 2003.  Revised on August 28, 2003}


\begin{abstract}
Let $E(1)_p$ denote the rational elliptic surface with a single
multiple fiber $f_p$ of multiplicity $p$.  We construct an
infinite family of homologous non-isotopic symplectic tori
representing the primitive $2$-dimensional homology class $[f_p]$
in $E(1)_p$ when\/ $p>1$. As a consequence, we get infinitely many non-isotopic
symplectic tori in the fiber class of the rational elliptic surface 
$E(1)\cong\cpk \#9\cpkk$. We also show how these tori can be
non-isotopically embedded as homologous symplectic submanifolds in
other symplectic $4$-manifolds.
\end{abstract}

\maketitle


\section{Introduction}

This paper is a continuation of studies initiated in \cite{ep1}
regarding infinite families of non-isotopic and symplectic tori
representing the same homology class in a symplectic 4-manifold.
Let $E(1)\cong\cpk \#9\cpkk$ be the rational elliptic surface
obtained by blowing up at nine branch points of a generic pencil of
cubic curves in $\cpk$.  Let $E(1)_p$ denote the complex surface
obtained by performing a logarithmic transformation of
multiplicity $p$\/ on a regular fiber $f$\/ of $E(1)$. Our main
result is the following:

\begin{theorem}\label{theorem:main}
For any integer\/ $p>1$, there exists an infinite family of
pairwise non-isotopic symplectic tori representing the primitive
homology class\/ $[f_p] \in H_2(E(1)_p)$, where\/ $f_p$ is the 
multiple fiber of $E(1)_p$.
\end{theorem}

Even though the tori we construct are not only non-isotopic but also 
inequivalent in the sense that there is no self-diffeomorphism of 
the ambient 4-manifold that carries one torus to another, we state 
the theorem the way we do since our interest in the subject is 
mainly due to the more general problem of finding non-isotopic symplectic 
representatives of a homology class, known as the symplectic isotopy
problem.

It is well-known that there exists a diffeomorphism between
$E(1)_p$ and $E(1)$, sending the homology class $[f_p]$ to $[f]\,$
(see e.g.$\;$\cite{gs}). Hence we immediately obtain the
following:

\begin{corollary}\label{cor:main}
For a suitable choice of a symplectic form on $E(1)$, there exists
an infinite family of pairwise non-isotopic symplectic tori
representing the primitive homology class\/ $[f] \in H_2(E(1))$,
where\/ $f$ is the the torus fiber of the elliptic fibration\/
$E(1)\rightarrow \mathbb{CP}^1$.
\end{corollary}


Examples of homologous, non-isotopic, symplectic tori were first
constructed in \cite{fs:non-isotopic} and then in \cite{ep1},
\cite{ep:k3} and \cite{vidussi:non-isotopic}. Infinite 
families of non-isotopic symplectic tori representing
divisible classes\/ $n[f_p]\in H_2(E(1)_p)$, $n\geq 2$, were
constructed in \cite{ep1}.  In \cite{ep:E(1)_K} and
\cite{vidussi:E(1)_K}, symplectic tori representing the same
primitive class in symplectic 4-manifolds homeomorphic but not
diffeomorphic to $E(1)$ were constructed.

It should be noted that the non-existence of such an infinite family 
of tori in the rational surfaces 
$\cpk$ and $\cpk \# \cpkk$ is proved by Sikorav in \cite{sikorav} 
and by Siebert and Tian in \cite{st}, respectively. Moreover, in contrast with
Corollary~\ref{cor:main}, which easily generalizes to 
$\cpk \# n \cpkk$ for $n \geq 9$, it is conjectured 
that there is at most one symplectic torus (up to isotopy) representing 
each homology class in $\cpk \# n \cpkk$ for $n < 9$.

The proof of Theorem~\ref{theorem:main} is spread out over the
next four sections. After reviewing the relevant definitions we interpret 
$E(1)_p$ as a link surgery manifold in
Section~\ref{sec:link surgery} and construct a family of homologous symplectic tori
in $E(1)_p$ in Section~\ref{sec:symplectic family}. In Section~\ref{sec:alexander}, 
Alexander polynomials of certain braid closures are given and then in 
Section~\ref{sec:sw} these polynomials are used to compute the Seiberg-Witten
invariants of fiber sums of $E(1)_p$ with $E(1)$, along the tori we construct
and a regular fiber $f$, respectively, which in turn give the non-isotopy of these tori. In
Section~\ref{sec:generalization}, we present a direct generalization 
in the form of Theorem~\ref{theorem:generalization}.  

In this introduction and elsewhere in the paper by isotopy we mean smooth
isotopy and all homology groups have integer coefficients.


\section{Link Surgery 4-Manifolds}
\label{sec:link surgery}

In this section, first we review the generalization of the
link surgery construction of
Fintushel and Stern \cite{fs:knots} by Vidussi
\cite{vidussi:smooth}, and then give specific link surgeries that
will be used in the following sections.

For an $n$-component link $L\subset S^3$,
choose an ordered homology basis of simple closed curves $\{(\alpha_i, \beta_i)
\}_{i=1}^{n}$ such that
the pair $(\alpha_i, \beta_i )$ lie in the $i$-th boundary
component of the link exterior
and the intersection of
$\alpha_i$ and $\beta_i$ is 1.
Let $X_i$
($i=1,\dots, n$) be a 4-manifold containing a 2-dimensional torus
submanifold $F_i$ of self-intersection $0$.  Choose a Cartesian product
decomposition $F_i = C_1^{i} \times C_2^{i}$, where each $C^i_j \cong S^1$
($j=1,2$) is an embedded
circle in $X_i$.

\begin{definition}\label{def:data}
The ordered collection
$$\mathfrak{D} \, =\; \big( \{(\alpha_i, \beta_i)
\}_{i=1}^{n}\: , \: \{(X_i  , \hspace{1pt}
F_i=  C_1^{i} \times C_2^{i} )\}_{i=1}^{n}\big)$$
is called a \emph{link
surgery gluing data}\/ for an $n$-component link $L$.
We define the \emph{link surgery manifold corresponding to} $\mathfrak{D}$
to be the closed $4$-manifold
\[
L(\mathfrak{D}) \: :=\; [\coprod_{i=1}^{n} X_i\setminus\nu
F_i]\hspace{-20pt}\bigcup_{F_i\times\partial D^2=(S^1\times
\alpha_i)\times\beta_i}\hspace{-20pt} [S^1\times(S^3\setminus \nu
L)]\, ,
\]
where $\nu$\/ denotes the tubular neighborhoods.  Here, the gluing
diffeomorphisms between the boundary 3-tori identify the torus
$F_i =  C_1^{i} \times C_2^{i}$\/ of $X_i$\/ with\/ $S^1\times
\alpha_i$\/ factor-wise, and act as the complex conjugation on the
last remaining\/ $S^1$ factor.
\end{definition}

\begin{remark}
Strictly speaking, the diffeomorphism type of the link surgery
manifold $L(\mathfrak{D})$\/ may possibly depend on the chosen
trivialization of $\,\nu F_i \cong F_i \times D^2$ (the framing of
$F_i$).  However, we will suppress this dependence in our
notation. It is well known (see e.g. \cite{gs}) that the
diffeomorphism type of $L(\mathfrak{D})$ is independent of the
framing of $F_i$ when $(X_i,F_i) = (E(1), f)$.
\end{remark}

We fix a Cartesian product decomposition of a regular torus
fiber\/ $f= C_1 \times C_2$ in the elliptic surface $E(1)$.  Also
fix\/ $T^2:=S^1\times S^1$.

\begin{figure}[!ht]
\begin{center}
\includegraphics[scale=.5]{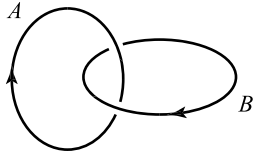}
\end{center}
\caption{Hopf link $L= A\cup B$}
\label{fig:hopf}
\end{figure}

\begin{definition}
Let\/ $L \subset S^3$ be the Hopf link in Figure~$\ref{fig:hopf}$.
For the link surgery gluing data\/
\begin{eqnarray}\label{eq:data}
 \mathfrak{D} \!\!\!  &:=& \!\!\!
\big(\{ (\mu(A),\lambda(A)), (\mu(B), \lambda(B)-p\mu(B)) \}, \\
&& \{ (X_1, F_1=C_1^1\times C_2^1), (T^2\times S^2, F_2 = T^2
\times \{ {\rm pt} \} )\}\big), \nonumber
\end{eqnarray}
we shall denote\/ $L(\mathfrak{D})$ by $(X_1)_p$.  Here,
$\mu(\,\cdot\,)$ and $\lambda(\,\cdot\,)$ denote the meridian and
the longitude of a knot, respectively.
\end{definition}

In particular, when $(X_1, F_1=C_1^1\times C_2^1) = (E(1),
f=C_1\times C_2)$, we denote $L(\mathfrak{D})$\/ by $E(1)_p$. This
notation is consistent with the existing literature as there is an
obvious diffeomorphism (see diffeomorphism (\ref{eq:cylinder}) in
the proof of Lemma~\ref{lemma:main}) between our link surgery
manifold\/
\begin{equation}\label{eq:decomposition}
L(\mathfrak{D})\:=\: [E(1)\setminus \nu
f]\,\cup\,[S^1\times(S^3\setminus L)]\,\cup\,[T^2\times D^2]\,
\end{equation}
and the logarithmic transform $E(1)_p = [E(1)\setminus \nu f]
\cup_{\varphi} [T^2\times D^2]$, where the gluing diffeomorphism
$\,\varphi: T^2\times
\partial D^2\rightarrow \partial(\nu f) $ induces the linear map
\begin{equation}\label{3x3 matrix}
\varphi_{\ast} \,= \left(\begin{array}{ccr}
  1\;\; & 0 & 0 \\
  0\;\; & 0 & 1 \\
  0\;\; & 1 & -p
\end{array}\right)
\end{equation}
between the first homology groups with the obvious choice of
bases.

\begin{lemma}\label{lemma:E(1)_p is symplectic}
If\/ $p\geq 1$, then $E(1)_p$ is a symplectic $4$-manifold
diffeomorphic to $E(1)$.  The diffeomorphism type of the
logarithmic transform\/ $E(1)_p$ does not depend on the choice of
the gluing map $\varphi$.  The homology class\/ $[f_p]\in
H_2(E(1)_p)$ is primitive.
\end{lemma}

\begin{proof}
See Theorem 3.3.3, Remark 3.3.5 and Theorem 8.3.11 in \cite{gs}.
\end{proof}


\section{Family of Homologous Symplectic Tori in $E(1)_p$}
\label{sec:symplectic family}

\begin{figure}[!ht]
\begin{center}
\includegraphics[scale=.85]{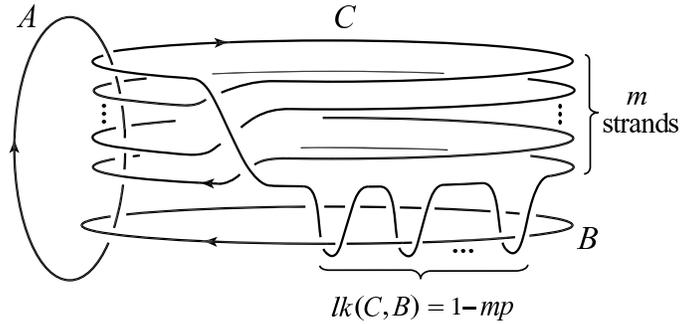}
\end{center}
\caption{3-component link $L_{m,p} =A\cup B \cup C\,$ in $S^3$}
\label{fig:3d}
\end{figure}

Let $T_C := S^1 \times C \subset [S^1 \times (S^3\setminus \nu L)]
\subset E(1)_p$, where the closed curve $C:=C_{m,p} \subset
(S^3\setminus \nu L )$\/ is given by Figure~\ref{fig:3d}. Note that
$C$ could be isotoped to an $(m,mp-1)$ torus knot in $\partial \nu (B) \subset S^3$.

\begin{lemma}\label{lemma:main}
For\/ every pair of integers\/ $p>1$\/ and $m\geq 1$, $T_C
=S^1\times C$ is a symplectic submanifold of\/ $E(1)_p$ and we
have\/ $[T_C]=[f_p]\,$ in\/ $H_2(E(1)_p)$.
\end{lemma}

\begin{proof}
To determine the homology class of $T_C$, note that\/ $[C]= m[\mu
(A)] + (1-mp)[\mu (B)]$\/ in\/ $H_1(S^3\setminus \nu L)$. When we glue\/
$(T^2\times D^2)$\/ to\/ $(S^1\times (S^3\setminus \nu L))$\/ using the gluing data
$\mathfrak{D}$ in (\ref{eq:data}), the homology class
$[\mu(A)]-p[\mu(B)]=[\lambda (B)]-p[\mu (B)]$ gets identified with $[\{{\rm
pt}\}\times\partial D^2] \in H_1(T^2\times D^2)$, which is
trivial.  Hence by K\"unneth's theorem, $[T_C]=[S^1\times
\mu(B)]\,$ in $H_2((S^1\times (S^3\setminus \nu L))\cup (T^2\times D^2))$.  Finally,
note that\/ $[S^1\times \mu(B)]=[f_p]\,$ and $[f]=[S^1\times
\mu(A)]=[S^1\times\lambda(B)]=[S^1 \times p\,\mu(B)]=p[f_p]\,$ in
$H_2(E(1)_p)$.

It is easy to see that the link exterior $S^3\setminus \nu L$\/ 
is diffeomorphic to $S^1 \times \mathbb{A}$, where
$\mathbb{A}\cong S^1 \times [0,1]\,$ is an annulus. Hence we may
identify the middle piece in the decomposition
(\ref{eq:decomposition}) above as the cylinder
\begin{equation}\label{eq:cylinder}
[S^1 \times (S^3\setminus \nu L)] \:\cong\:
[S^1 \times (S^1 \times \mathbb{A})]
\:\cong\: T^3\times [0,1] .
\end{equation}
Since we assume that $p>1$, the logarithmic transform corresponds to the 
rational blowdown construction of Fintushel and Stern \cite{fs:blowdown}. On the
other hand, by \cite{symington}, rational blowdown can be done
symplectically and the symplectic form needs to be modified only near 
the collar neighborhood of the boundary lens spaces involved. As a result of this
and since $f \simeq p\cdot f_p$ remains symplectic in $E(1)_p$,
we may assume that the symplectic form on $E(1)_p$ restricts to
\[
\omega \,=\, dx\wedge dy \,+\, r\,dr\wedge d\theta
\]
on the subset $\,T^3\times[\epsilon,1]\subset[S^1 \times
(S^1\times \mathbb{A})]$ for suitable\/ $0<\epsilon< 1$, where
$x$\/ and $y$\/ are the angular coordinates on the first and the
second $S^1$ factors respectively, and $(r,\theta)$ are the polar
coordinates on the annulus $\mathbb{A}$.

We can embed the curve $C$\/ inside
$\,T^2\times[\epsilon,1]\subset(S^1\times \mathbb{A})$ such that
$C$\/ is transverse to every annulus of the form, $\{{\rm pt}\}
\times \mathbb{A}$, and the restriction $dy|_{C}$ never vanishes.
It follows that $\omega |_{T_C} = (dx\wedge dy)|_{T_C} \neq 0$,
and consequently $T_C$ is a symplectic submanifold of $E(1)_p$.
\end{proof}


\section{Alexander Polynomials of Certain Braid Closures}
\label{sec:alexander}

\begin{lemma}\label{lem:alex}
$\Delta_{L_{m,p}}(x,s,t)=$ 
$$1-x(st)^{1-mp}+xt\frac{(xt)^{m-1}-1}{xt-1}
\left[1+(s-1)\frac{(st)^{1-mp}-1}{st-1}-x(st)^{1-mp}\right] \, , $$ 
where $L_{m,p}$ is the 3-component link in Figure~\ref{fig:3d}, and
the variables\/ $x$, $s$ and\/
$t$\/ correspond to the components\/ $A$, $B$\/ and\/ $C$\/
respectively.
\end{lemma}

\begin{figure}[!ht]
\begin{center}
\includegraphics[scale=.65]{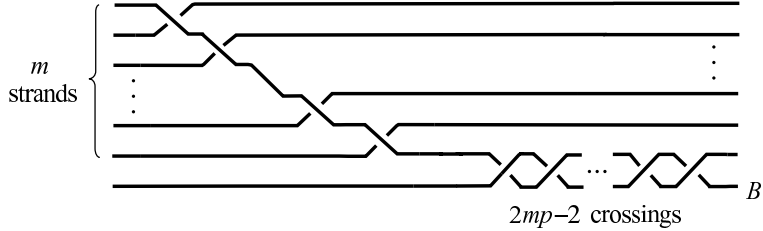}
\end{center}
\caption{$(m+1)$-strand braid whose closure gives $B\cup C$}
\label{fig:torusknotbraid}
\end{figure}

\begin{proof}
This follows from the formula in Theorem 1 of
\cite{morton} which gives the multivariable Alexander polynomial
of a braid closure and its axis in terms of the representation of
the braid. We view $A$\/ as the axis of the $(m+1)$-strand braid
in Figure~\ref{fig:torusknotbraid}. Note that the closure of this braid
and its axis is $L_{m,p}=A\cup B\cup C$. This statement can also be thought of
as a slightly modified special case of Lemma 6.2 in \cite{ep:k3}. 
\end{proof}

\begin{remark}\label{rem:symalex}
Note that the symmetrized Alexander polynomial of $L_{m,p}$ is
$$\Delta^{\rm sym}_{L_{m,p}} (x,s,t) = x^{-m/2}s^{(mp-1)/2}t^{(mp-m)/2} 
\cdot \Delta_{L_{m,p}} (x,s,t)\, .$$
\end{remark}


\section{Non-Isotopy:  Seiberg-Witten Invariants}
\label{sec:sw}

Our strategy is to show that for fixed choice of $p$ the isotopy
types of the tori $\{ T_C \}_{m\geq 1}$ can be distinguished by
comparing the Seiberg-Witten invariants of the corresponding
family of fiber sum 4-manifolds $\{E(1)_p \#_{T_C=f} E(1)\}_{m\geq
1}$. Note that there is a canonical framing of a regular fiber
$f$\/ in $E(1)$, coming from the elliptic fibration
$E(1)\rightarrow \mathbb{CP}^1$.

\begin{lemma}
The fiber sum\/ $E(1)_p \#_{T_C=f} E(1)$\/ is diffeomorphic to the
link surgery manifold $L_{m,p}(\mathfrak{D}')$, where
\begin{eqnarray}
&&  \quad\quad   \mathfrak{D}' : = \:  \big(
\{(\mu(A),\lambda(A)),(\mu(B),\lambda(B)-p\hspace{1pt}\mu(B)),
(\lambda(C),-\mu(C))\},  \\[3pt]
&& \{(E(1),f=C_1\times C_2), (T^2\times S^2,F_2=S^1 \times S^1
\times\{{\rm pt} \}),(E(1), f =C_1\times C_2)\}\big). \nonumber
\end{eqnarray}
\end{lemma}

\begin{proof}
It was already observed in \cite{ep:k3} that the fiber sum
construction corresponds to this type of link surgery.  (See also
\cite{ep:E(1)_K}.)
\end{proof}

Recall that the Seiberg-Witten invariant\/
$\overline{SW}_{\!\!X}$\/ of a 4-manifold $X$\/ can be thought of
as an element of the group ring of $H_2(X)$, i.e.\/
$\overline{SW}_{\!\!X} \in \zz [ H_2( X ) ]$. If we write\/
$\overline{SW}_{\!\!X} = \sum_g a_g g \hspace{1pt}$, then we say
that\/ $g\in H_2( X )$\/ is a Seiberg-Witten \emph{basic class}\/
of $X$\/ if\/ $a_g\neq 0$. Since the Seiberg-Witten invariant of a
4-manifold is a diffeomorphism invariant, so are the
divisibilities of Seiberg-Witten basic classes. The Seiberg-Witten
invariant of the link surgery manifold $L_{m,p}(\mathfrak{D}')$ is
known to be related to the Alexander polynomial\/
$\Delta_{L_{m,p}}$ of the link $L_{m,p}$.

\begin{theorem}\label{theorem:seiberg-witten}
Let\/ $\iota:[S^1\times(S^3\setminus\nu L_{m,p})]\rightarrow
L_{m,p}(\mathfrak{D}')$\/ be the inclusion map.  Let\/
$\varsigma:=\iota_{\ast}[S^1\times\mu(B)],$
$\tau:=\iota_{\ast}[S^1\times\mu(C)] \in
H_2(L_{m,p}(\mathfrak{D}') ).$  Then\/ $\varsigma$ and\/ $\tau$
are 
linearly independent. The Seiberg-Witten
invariant of\/ $L_{m,p}(\mathfrak{D}')$ is given by


\begin{eqnarray}\label{eq:sw}
\overline{SW}_{\!L_{m,p}(\mathfrak{D}')} &=& \varsigma^{-mp+2+2p}\tau^{-mp^2+mp} \cdot \\
&& 
\Big[ 
( \varsigma^{2mp-2p-4} + \cdots + \varsigma^2 + 1 ) \cdot
\nonumber \\
&& 
(\varsigma^{2mp-2p}\tau^{2m^2p-2mp} + \cdots + \varsigma^{2p}\tau^{2mp} + 1 ) \nonumber \\
&&
+\tau^2 ( \varsigma^{2mp-4}\tau^{2mp-4} + \cdots + \varsigma^2\tau^2 + 1  ) \cdot
\nonumber \\
&&  
( \varsigma^{2mp-4p}\tau^{2m^2p-4mp} + \cdots + \varsigma^{2p}\tau^{2mp} + 1  ) \Big]
 \nonumber \, \, . 
\end{eqnarray}

\end{theorem}

\begin{proof}
Let\/ $N := (S^3\setminus \nu L_{m,p})$, and let $Z:=(T^2\times
D^2)$. Recall from \cite{doug:pft3} that we have
$\overline{SW}_{\!E(1)\setminus \nu f}=1$, and also
$$
\overline{SW}^{\,\pm}_{\!(T^2 \times D^2)} =\:
\frac{1}{[T^2\times\{ {\rm pt} \}]^{-1}-[T^2\times\{{\rm pt} \}]}
\; .$$

From the gluing formulas in \cite{doug:pft3} and
\cite{Taubes:T^3}, we may conclude that
\begin{eqnarray}\label{eqn:prod}
\overline{SW}_{\!L_{m,p}(\mathfrak{D}')} &=&
\overline{SW}_{\!E(1)\setminus \nu f} \cdot
\overline{SW}_{\!E(1)\setminus \nu f} \cdot
\overline{SW}^{\,\pm}_{\!(T^2\times D^2)} \cdot
\Delta_{L_{m,p}}^{\rm sym}(\xi^2,\varsigma^2,\tau^2) \nonumber \\
&=&  \frac{1}{[T^2\times\{ {\rm pt} \}]^{-1}-[T^2\times\{{\rm pt} \}]} \cdot
\Delta_{L_{m,p}}^{\rm sym}(\xi^2,\varsigma^2,\tau^2) \, , 
\end{eqnarray}
where $\xi:=\iota_{\ast}[S^1\times\mu(A)]$.

First of all, $[\lambda (B)] =  [\mu (A)] + (1-mp) [\mu (C)] \in H_1(N)$. On the
other hand, according to the gluing data $\mathfrak{D}'$, 
$\lambda (B) - p \mu (B)$ is identified with $\{ {\rm pt} \} \times \partial D^2 \in Z$, therefore
$[\lambda (B)] = p [\mu (B)]$ in $L_{m,p}(\mathfrak{D}')$
and consequently 
\begin{equation}\label{eqn:xi}
\xi = \varsigma^p\tau^{mp-1} \, .
\end{equation} 

Also note that 
$T^2 \in T^2 \times S^2$ and $S^1 \times \mu (B) $ 
are identified by the gluing data $\mathfrak{D}'$, hence
\begin{equation}\label{eqn:t2}
[T^2 \times \{ pt \} ] = \varsigma \in H_2 (L_{m,p}(\mathfrak{D}')) \, .
\end{equation} 

From Lemma~\ref{lem:alex} and Remark~\ref{rem:symalex} we have 
\begin{eqnarray}\label{eqn:symalex}
\Delta_{L_{m,p}}^{\rm sym}((\varsigma^{p}\tau^{mp-1})^2,\varsigma^2,\tau^2) &=& 
\varsigma^{-1}\tau^{-m^2p+mp} \cdot
\bigg[ 1-\varsigma^{-2mp+2p+2} \\
\hspace{1cm} 
+\varsigma^{2p}\tau^{2mp} \frac{(\varsigma^{2p}\tau^{2mp})^{m-1}-1}{\varsigma^{2p}\tau^{2mp}-1} 
&& \hspace{-1cm} \left(1+(\varsigma^2-1)\frac{(\varsigma^2\tau^2)^{1-mp}-1}{\varsigma^2\tau^2-1}
-\varsigma^{-2mp+2p+2}\right)  \bigg] \, . \nonumber
\end{eqnarray}

After putting Equations (\ref{eqn:prod}), (\ref{eqn:xi}), 
(\ref{eqn:t2}) and (\ref{eqn:symalex}) together and simplifying 
the resulting expression, we get (\ref{eq:sw}). 

Next we show that $\varsigma$ and $\tau$ are 
linearly independent elements of $H_2(L_{m,p}(\mathfrak{D}'))$. A
Mayer-Vietoris argument, combined with Freedman's classification
theorem (see, for example, $\;$\cite{FQ}), shows that $L_{m,p}(\mathfrak{D}')$\/
is homeomorphic to $E(2)$.

Consider the composition of homomorphisms
\begin{equation}\label{eq:homomorphism}
H_1(N) \longrightarrow H_2(S^1\times N)
\stackrel{\iota_{\ast}\;}{\longrightarrow}
H_2(L_{m,p}(\mathfrak{D}')),
\end{equation}
where the first
map is a part of the K\"unneth isomorphism
\begin{equation}\label{eq:Kunneth}
H_1(N) \oplus H_2(N)
\stackrel{\cong}{\longrightarrow} H_2(S^1\times N).
\end{equation}
Note that $H_2(N) \cong \zz\oplus \zz\,$,
as is easily seen from the long exact
sequence of the pair\/ $(N,\partial N)$ as follows.
\[
\begin{array}{ccccccccc}
H^0(N)&\longrightarrow & H^0(\partial N) &\longrightarrow &
 H^1(N,\partial N) & \stackrel{0}{\longrightarrow} &
H^1(N) & \longrightarrow &  H^1(\partial N) \\
_{||} && _{||} && _{||} & & _{||}
&& _{||}\\[3pt]
\zz  & \longrightarrow  & \zz^3 & \longrightarrow & H_2(N) &
 \stackrel{0}{\longrightarrow}
& \zz^3 &
\longrightarrow & \zz^6
\end{array}
\]
Note that the first map sends the generator $\hspace{1pt} 1\in \zz\,$ to
the diagonal element\/ $(1,1,1)\in \zz^3$, while the last map is injective.
We have also used the Lefschetz duality theorem (for manifolds with boundary)
to identify $H_2(N)\cong
H^1(N,\partial N)$.

Next consider the long exact sequence of the pair
$(L_{m,p}(\mathfrak{D}'), S^1\times N)$:
\[
0 = H_3(L_{m,p}(\mathfrak{D}')) \longrightarrow
H_3(L_{m,p}(\mathfrak{D}'),S^1\times N) \longrightarrow
H_2(S^1\times N) \stackrel{\iota_{\ast}\;}{\longrightarrow}
H_2(L_{m,p}(\mathfrak{D}'))
\]
The kernel of the last map $\iota_{\ast}$ is isomorphic to
$H_3(L_{m,p}(\mathfrak{D}'),S^1\times N)$.  By Lefschetz duality
theorem (for relative manifolds),
$H_3(L_{m,p}(\mathfrak{D}'),S^1\times N)$ is in turn isomorphic to
\[
H^1(L_{m,p}(\mathfrak{D}')\setminus (S^1\times N)) \:\cong\;
H^1(E(1)\setminus \nu f) \oplus H^1(E(1)\setminus \nu f) \oplus
H^1(Z).
\]
Since we have $H^1(E(1)\setminus \nu f) =0\,$ and
\begin{equation}\label{eq:z plus z}
H^1(Z) =H^1(T^2\times D^2 ) \cong \; \zz\oplus \zz \, ,
\end{equation}
the kernel of $\iota_{\ast}$ is isomorphic to\/ $\zz\oplus\zz\,$.

Finally we observe that only one $\zz$ summand of (\ref{eq:z plus z}) lies
in the kernel of the composition (\ref{eq:homomorphism}).  The other
$\zz$ summand belongs to the kernel of
\[
H_2(N) \longrightarrow H_2(S^1\times N)
\stackrel{\iota_{\ast}\;}{\longrightarrow}
H_2(L_{m,p}(\mathfrak{D}')),
\]
where the first map is the second part of the K\"unneth
isomorphism (\ref{eq:Kunneth}). We have thus shown that the kernel
of the composition (\ref{eq:homomorphism}) is of rank one.  It
follows immediately that $\varsigma$ and $\tau$ are linearly
independent, since we have already showed that $\xi=\varsigma^p\tau^{mp-1}$. 
A more detailed analysis shows that\/ $\{\varsigma ,\tau\}$\/ can be
extended to a basis of $H_2(L_{m,p}(\mathfrak{D}'))$, which we
shall omit. (Also see the proof of Proposition 3.2 in
\cite{McMullen-Taubes} for a similar argument.)
\end{proof}

\begin{corollary}\label{cor:basics}
The number of Seiberg-Witten basic classes of $L_{m,p}(\mathfrak D)$ is
$$ 2pm^2-2(p+1)m+1 \, .$$
\end{corollary}

\begin{proof}
First of all, the summands in Equation(\ref{eq:sw}) have no common terms, i.e.
terms with matching powers for $\varsigma$ and $\tau$, because the general term 
in the first summand is $(\varsigma^2)^k(\varsigma^{2p}\tau^{2mp})^l$, where 
$0 \leq k \leq mp-p-2$ and $0 \leq l \leq m-1$, the general term in the second 
summand is $\tau^2(\varsigma^2\tau^2)^i(\varsigma^{2p}\tau^{2mp})^j$, where
$0 \leq i \leq mp-2$ and $0 \leq j \leq m-2$, and in order for two such terms 
to have matching powers, for even only $\tau$, $2lmp$ should be the same as
$2+2i+2jmp$ which would imply $i \equiv -1 \, (\mbox{mod} \, mp)$ contradicting
$0 \leq i \leq mp-2$. Therefore, to get the number of basic classes, we simply add
the number of terms in each of these summands.
\end{proof}

\begin{corollary}\label{cor:non-isotopy}
There is no self-diffeomorphism of\/ $E(1)_p$ that maps one element 
of the family of tori\/ $\{ T_C \}_{m \geq 1}$ to another. In particular,
these tori are pairwise non-isotopic.
\end{corollary}

\begin{proof}
Since the Seiberg-Witten
invariant is a diffeomorphism invariant, so are the number
of basic classes. On the other hand, by Corollary~\ref{cor:basics}, 
the number of basic classes of $L_{m,p}(\mathfrak{D}')$ is 
$ 2pm^2-2(p+1)m+1$ and for fixed $p\geq 2$ this polynomial of $m$ is 
obviously increasing when $m \geq 1$. Hence $L_{m,p}(\mathfrak{D}')$ is
diffeomorphic to $L_{m',p}(\mathfrak{D}')$ if and only if\/ $m=m'$
proving that the tori in $\{ T_C \}_{m \geq 1}$ are different up
to isotopy and in fact even up to self-diffeomorphisms of\/
$E(1)_p$.
\end{proof}

This concludes the proof of Theorem~\ref{theorem:main}.


\section{Generalization to Other Symplectic 4-Manifolds}
\label{sec:generalization}

We can easily extend Theorem~\ref{theorem:main} to $E(n)_p$
$(n\geq 2)$ and more generally to $X_p$, where $X$\/ is a
symplectic 4-manifold satisfying certain topological conditions as
in \cite{ep1}.

\begin{theorem}\label{theorem:generalization}
Assume that\/ $F$ is a symplectic\/ $2$-torus in a symplectic\/
$4$-manifold\/ $X$. Suppose that\/ $[F]\in H_2(X)$ is primitive,
$[F]\cdot[F]=0$, and that $F$\/ lies in a fishtail neighborhood.
If $\,b_2^+(X)=1$, then we also assume that $\,\overline{SW} _{\!
X\setminus \nu F}\neq 0\,$ and is a finite sum. Then there exists
an infinite family of pairwise non-isotopic symplectic tori in\/
$X_p$ representing the homology class\/ $[F_p]\in H_2(X_p)$ for
any integer $p>1$.  These tori are inequivalent under
self-diffeomorphisms of $X_p$.
\end{theorem}

\begin{proof}
It was shown in \cite{fs:blowdown} and \cite{symington} that $X_p$
possesses a canonical symplectic form coming from the symplectic
form on $X$.  The rest of the proof goes the same way as before.
\end{proof}


\subsection*{Acknowledgments}
We would like to thank Ronald Fintushel and Stefano Vidussi for
their encouragement and helpful comments on this and other works
of ours.


\end{document}